\documentclass[12pt]{amsart}
\usepackage{amsmath, amsthm, amscd,
amsfonts}

\newtheorem{theorem}{Theorem}[section]
\newtheorem{lemma}[theorem]{Lemma}
\newtheorem{proposition}[theorem]{Proposition}
\newtheorem{corollary}[theorem]{Corollary}
\theoremstyle{definition}
\newtheorem{definition}[theorem]{Definition}
\newtheorem{example}[theorem]{Example}

\theoremstyle{remark}

\numberwithin{equation}{section}
\newfont{\aj}{eufm10 at12pt}
\newfont{\ajk}{eufm10 at10pt}
\newcommand{\cala}{\mbox{\aj A}}
\newcommand{\calb}{\mbox{\aj B}}
\newcommand{\calh}{\mbox{\aj H}}
\newcommand{\calx}{\mbox{\aj X}}
\newcommand{\call}{\mbox{\aj L}}
\newcommand{\calk}{\mbox{\aj K}}
\newcommand{\cals}{\mbox{\aj S}}
\newcommand{\calak}{\mbox{\ajk A}}
\newcommand{\calhk}{\mbox{\ajk H}}
\begin{document}
\title[Achievement of continuity of
$(\varphi,\psi)$-derivations
]{ Achievement of continuity of
$(\varphi,\psi)$-derivations without linearity}
\author[S. Hejazian, A.R. Janfada, M. Mirzavaziri, M.S.
Moslehian]{S. Hejazian, A. R. Janfada, M. Mirzavaziri and M. S.
Moslehian}
\address{Shirin Hejazian: Department of Mathematics,
Ferdowsi University, P. O. Box 1159, Mashhad 91775, Iran; \newline
Banach Mathematical Research Group (BMRG), Mashhad, Iran.}
\email{hejazian@ferdowsi.um.ac.ir}
\address{Ali Reza Janfada: Department of Mathematics, Ferdowsi University, P. O. Box 1159,
Mashhad 91775, Iran} \email{ajanfada@math.um.ac.ir}
\address{Madjid Mirzavaziri: Department of Mathematics, Ferdowsi University, P. O.
Box 1159, Mashhad 91775, Iran;\newline Centre of Excellence in
Analysis on Algebraic Structures (CEAAS), Ferdowsi Univ., Iran.}
\email{mirzavaziri@math.um.ac.ir}
\address{Mohammad Sal Moslehian: Department of Mathematics, Ferdowsi University, \newline  P. O.
Box 1159, Mashhad 91775, Iran;\newline Centre of Excellence in
Analysis on Algebraic Structures (CEAAS), Ferdowsi Univ., Iran.}
\email{moslehian@ferdowsi.um.ac.ir} 

\begin{abstract}Suppose that $\calak$ is a $C^*$-algebra acting
on a Hilbert space $\calhk$, and that $\varphi, \psi$ are mappings
from $\calak$ into $B(\calhk)$ which are not assumed to be
necessarily linear or continuous. A $(\varphi, \psi)$-derivation
is a linear mapping $d: \calak \to B(\calhk)$ such that
$$d(ab)=\varphi(a)d(b)+d(a)\psi(b)\quad (a,b\in \calak).$$ We prove that if $\varphi$ is a multiplicative (not
necessarily linear)\ $*$-mapping, then every
$*$-$(\varphi,\varphi)$-derivation is automatically continuous.
Using this fact, we show that every
$*$-$(\varphi,\psi)$-derivation $d$ from $\calak$ into
$B(\calhk)$ is continuous if and only if the $*$-mappings
$\varphi$ and $\psi$ are left and right $d$-continuous,
respectively.
\end{abstract}
\maketitle

\section{Introduction}

Recently, a number of analysts \cite{B-M, B-V, M-M, M-M1, MOS}
have studied various generalized notions of derivations in the
context of Banach algebras. There are some applications in the
other fields of research \cite{H-L-S}. Such mappings have been
extensively studied in pure algebra; cf. \cite{A-R, BRE, HVA}. A
generalized concept of derivation is as follows.

\begin{definition}\label{maindef}
Suppose that $\calb$ is an algebra, $\cala$ is a subalgebra of
$\calb$, $\calx$ is a $\calb$-bimodule, and $\varphi,\psi :\cala
\to \calb$ are mappings. A linear mapping $d:\cala \to \calx$ is a
\emph{$(\varphi,\psi)$-derivation} if
$$d(ab)=\varphi(a)d(b)+d(a)\psi(b)\quad (a,b\in \cala).$$ By a
\emph{ $\varphi$-derivation} we mean a
$(\varphi,\varphi)$-derivation. Note that we do not have any extra
assumptions such as linearity or continuity on the mappings
$\varphi$ and $\psi$.
\end{definition}
The automatic continuity theory is the study of (algebraic)
conditions on a category, e.g. $C^*$-algebras, which guarantee
that every mapping belonging to a certain class, e.g.
derivations, is continuous. S. Sakai \cite{SAK} proved that every
derivation on a $C^*$-algebra is automatically continuous. This
is an affirmative answer to the conjecture made by I. Kaplansky
\cite{KAP}. J. R. Ringrose \cite{RIN} extended this result for
derivations from a $C^*$-algebra $\cala$ to a Banach
$\cala$-bimodule $\calx$. B. E. Johnson and A. M. Sinclair
\cite{J-S} proved that every derivation on a semisimple Banach
algebra is automatically continuous. Automatic continuity of
module derivations on $JB^*$-algebras have been studied in
\cite{H-N}. The reader may find a collection of results
concerning these subjects in \cite{DAL2, SIN, VIL}.

M. Bre\v{s}ar and A. R. Villena \cite{B-V} proved that for an
inner automorphism $\varphi$, every $(id, \varphi)$-derivation on
a semisimple Banach algebra is continuous, where $id$ denotes the
identity mapping. In \cite{M-M}, it is shown that every
$(\varphi,\psi)$-derivation from a $C^*$-algebra $\cala$ acting
on a Hilbert space $\calh$ into $B(\calh)$ is automatically
continuous, if $\varphi$ and $\psi$ are continuous $*$-linear
mappings from $\cala$ into $B(\calh)$.

This paper consists of five sections. We define the notion of a
$(\varphi,\psi)$-derivation in the first section and give some
examples in the second. In the third section, by using methods of
\cite{M-M}, we prove that if $\varphi$ is a multiplicative (not
necessarily linear) $*$-mapping from a $C^*$-algebra $\cala$
acting on a Hilbert space $\calh$ into $B(\calh)$, then every
$*$-$\varphi$-derivation $d: \cala \to B(\calh)$ is automatically
continuous. In the fourth section, we show that for not
necessarily linear or continuous $*$-mappings $\varphi,
\psi:\cala \to B(\calh)$ the continuity of a
$*$-$(\varphi,\psi)$-derivation $d: \cala \to B(\calh)$ is
equivalent to the left $d$-continuity of $\varphi$ and the right
$d$-continuity of $\psi$. The mapping $\varphi$ (resp. $\psi$) is
called left (resp. right) $d$-continuous if
$\displaystyle{\lim_{x \to 0}} \big(\varphi(x)d(b)\big)=0$, for
all $b\in \cala$ (resp. $\displaystyle{\lim_{x \to
0}}\big(d(b)\psi(x)\big)=0$, for all $b\in \cala$). Obviously
these conditions happen whenever $\displaystyle{\lim_{x\to
0}}\varphi(x)=0=\displaystyle{\lim_{x\to 0}} \psi(x)$, in
particular whenever $\varphi$ and $\psi$ are bounded linear
mappings. Thus we extend the main results of \cite{M-M} to a
general framework. Furthermore, we prove that if $d$ is a
continuous $*$-$(\varphi,\psi)$-derivation, then we can replace
$\varphi$ and $\psi$ with mappings with `at most' zero separating
spaces. The last section is devoted to study the continuity of the
so-called generalized $*$-$(\varphi,\psi)$-derivations from
$\cala$ into $B(\calh)$.

The reader is referred to \cite{DAL2} for undefined notations and
terminology.

\section{Examples}

In this section, let $\calb$ be an algebra, let $\cala$ be a
subalgebra of $\calb$, and let $\calx$ be a $\calb$-bimodule. The
following are some examples concerning Definition \ref{maindef}.

\begin{example}
Every ordinary derivation $d: \cala \to \calx$ is an
$id$-derivation, where $id: \cala \to \cala$ is the identity
mapping.
\end{example}

\begin{example}
Every homomorphism $\rho :\cala \to\cala$ is a
$(\frac{\rho}{2},\frac{\rho}{2})$-derivation.\end{example}

\begin{example}Let
$\varphi,\psi:\cala\to\cala$ be homomorphisms, and let $x\in\calx$
be a fixed element. Then the linear mapping $d_x:\cala\to\calx$
 defined by  $$d_x(a) :=\varphi(a)x - x\psi(a)
\qquad(a\in \cala),$$ is a $(\varphi,\psi)$-derivation, which is
called an inner $(\varphi,\psi)$-derivation corresponding to $x$.
\end{example}

\begin{example}\label{example}
Assume that  $\gamma, \theta:\mathcal{C}[0,2]\to \mathcal{C}[0,2]$
are arbitrary mappings, where $\mathcal{C}[0,2]$ denotes the
$C^*$-algebra of all continuous complex valued functions on
$[0,2]$. Take $\lambda \in\mathbb{C}$ and fixed elements
$f_1,f_2$,  and $h_0$ in $\mathcal{C}[0,2]$, such that $$f_1
h_0=0=f_2 h_0.$$ For example, let us take \begin{eqnarray*}h_0
(t)&:=&(1-t)\chi_{_{[0,1]}}(t),\\
f_1(t)&:=&(t-1)\chi_{_{[1,2]}}(t),\\
f_2(t)&:=&(t-\frac{3}{2})\chi_{_{[\frac{3}{2},2]}}(t),\end{eqnarray*}
where $\chi_{_E}$ denotes the characteristic function of $E$.
Define $\varphi,\psi,d:\mathcal{C}[0,2]\to\mathcal{C}[0,2]$ by
\begin{eqnarray*}
\varphi(f)&:=&\lambda
f+\gamma(f)f_1,\\\psi(f)&:=&(1-\lambda)f+\theta(f)f_2,\\ d(f)&:=&f
h_0.
\end{eqnarray*}
Then $d$ is a $(\varphi,\psi)$-derivation, since
\begin{eqnarray*}
\varphi(f)d(g)+d(f)\psi(g)&=&\big(\lambda f+\gamma(f)f_1\big)(g
h_0)+(fh_0)\big((1-\lambda)g+\theta(f)f_2\big)\\ &= &\lambda
fgh_0 + \gamma(f)gf_1h_0 + (1-\lambda)fgh_0 + \theta(f)fh_0f_2\\
&= &fgh_0\\ &= &d(fg).
\end{eqnarray*}
Moreover, we may choose $\gamma$ and $\theta$ such that $\varphi$
and $\psi$ are neither linear nor continuous.
\end{example}

\section{multiplicative mappings}

In this section, we are going to show  how a multiplicative
property gives us the linearity. We start our work with some
elementary properties of $(\varphi,\psi)$-derivations.

\begin{lemma}\label{algebraic}
Let $\calb$ be an algebra, let $\cala$ be a subalgebra of
$\calb$, and let $\calx$ be a $\calb$-bimodule. If $d:\cala \to
\calx$ is a $(\varphi,\psi)$-derivation,
then\\$(i)~~~\big(\varphi(ab)-\varphi(a)\varphi(b)\big)d(c)=d(a)\big(\psi(bc)-\psi(b)\psi(c)\big)$;\\
$(ii)~~~\big(\varphi(a+b)-\varphi(a)-\varphi(b)\big)d(c)=0=d(a)\big(\psi(b+c)-\psi(b)-\psi(c)\big)$;\\
$(iii)~~~\big(\varphi(\lambda
a)-\lambda\varphi(a)\big)d(b)=0=d(a)\big(\psi(\lambda
b)-\lambda(\psi(b)\big)$,\\for all $a,b,c\in\cala$, and all
$\lambda \in \mathbb{C}$.
\end{lemma}

\begin{proof}
Let $a,b,c\in\cala$ and $\lambda \in \mathbb{C}$. For the first
equation we have
\begin{eqnarray*}
0&=&d\big((ab)c\big)-d\big(a(bc)\big)\\&=&\varphi(ab)d(c)+d(ab)\psi(c)-\varphi(a)d(bc)-d(a)\psi(bc)\\
&=&\varphi(ab)d(c)+\big(\varphi(a)d(b)+d(a)\psi(b)\big)\psi(c)\\
&&-\varphi(a)\big(\varphi(b)d(c)+d(b)\psi(c)\big)-d(a)\psi(bc)\\
&=&\big(\varphi(ab)-\varphi(a)\varphi(b)\big)d(c)-d(a)\big(\psi(bc)-\psi(b)\psi(c)\big).
\end{eqnarray*}
Calculating $d\big((a+b)c\big)-d(bc)-d(ac)$ and $d\big(a(\lambda
b)\big)-\lambda d(ab)$, we obtain the other equations.
\end{proof}

We also need the following lemma. Recall that, for a Hilbert
space $\calh$, a subset $Y$ of $B(\calh)$ is said to be
self-adjoint if $u^*\in Y$ for each $u\in Y$.

\begin{lemma}\label{L perp=}
Let $\calh$ be a Hilbert space, and let $Y$ be a self-adjoint
subset of $B(\calh)$. Assume that
$\call_0=\displaystyle{\bigcup_{u\in Y}}u(\calh)$, and $\call$ is
the closed linear span of $\call_0$. If $\calk=\call^\perp$, then
$\calk=\displaystyle{\bigcap_{u\in Y}}\ker(u)$.
\end{lemma}

\begin{proof}
Let $k\in \calk$. Then $\langle \ell,k\rangle=0$, for all
$\ell\in\call$. Since $Y$ is self-adjoint, $$\langle u(k),h\rangle
=\langle k,u^*(h)\rangle =0$$ for all $h\in \calh$ and all $u\in
Y$. Therefore $u(k)=0$, for all $u\in Y$. Hence
$\calk\subseteq\displaystyle{\bigcap_{u\in Y}}\ker(u)$. The
inverse inclusion can be proved similarly.
\end{proof}

Recall that a mapping $f:\cala \to B(\calh)$ is called a
$*$-mapping, if $f(a^*)=f(a)^*$, for all $a\in \cala$. If we
define $f^*:\cala \to B(\calh)$, by $f^*(a)
:=\big(f(a^*)\big)^*$, then $f$ is a $*$-mapping if and only if
$f^*=f$.

\begin{theorem}\label{multipl}
Suppose that $\cala$ is a $C^*$-algebra acting on a Hilbert space
$\calh$. Let $\varphi:\cala\to B(\calh)$ be a $*$-mapping and let
$d: \cala\to B(\calh)$ be a $*$-$\varphi$-derivation. If $\varphi$
is multiplicative, i.e. $\varphi(ab) = \varphi(a)\varphi(b)
\quad(a, b \in \cala)$, then $d$ is continuous.
\end{theorem}

\begin{proof}
Set
\begin{eqnarray*}Y :=
\{\varphi(\lambda b)-\lambda\varphi(b)\mid \lambda \in {\mathbb
C}, b\in \cala\}\cup \{\varphi(b+c)-\varphi(b)-\varphi(c)\mid
b,c\in \cala\}.
\end{eqnarray*}
Since $\varphi$ is a $*$-mapping, $Y$ is a self-adjoint subset of
$B(\calh)$. Put
$$\call_0 :=\bigcup\{u(h)\mid u\in Y, h\in \calh\}.$$ Let $\call$
be the closed linear span of $\call_0$, and let $\calk
:=\call^\perp$. Then $\calh=\calk\oplus\call$, and by Lemma \ref{L
perp=}, we get\begin{eqnarray}\label{ker'}\calk=\bigcap
\{\ker(u)\mid u\in Y\}.\end{eqnarray}Suppose that $p\in B(\calh)$
is the orthogonal projection of $\calh$ onto $\calk$, then by
(\ref{ker'}) we have
\begin{eqnarray}\label{scal'}\varphi(\lambda
b) p&=&\lambda\varphi(b) p,
\end{eqnarray}
\begin{eqnarray}\label{add'}\varphi(b+c) p =
\varphi(b) p+\varphi(c) p,
\end{eqnarray}
for all $b,c\in \cala$, and all $\lambda\in \mathbb{C}$. Now we
claim that,
\begin{eqnarray}\label{commut'}
p\varphi(a)=\varphi(a)p\,,\,\,\,pd(a)=d(a)=d(a)p\qquad (a\in
\cala).\end{eqnarray} To end this, note that by Lemma
\ref{algebraic}, we
have\begin{eqnarray*}&&d(a)\big(\varphi(\lambda b)
-\lambda\varphi(b)\big)=0,\\&&d(a)\big(\varphi(b+c)
-\varphi(b)-\varphi(c)\big)=0,\end{eqnarray*} for all $a, b, c \in
\cala$, and all $\lambda\in \mathbb{C}$. Thus $d(a)\call=0$, and
so $d(a)(1-p)=0$, for all $a \in \cala$. Since $d$ is a
$*$-mapping, it follows that $pd(a)=d(a)=d(a)p$, for all $a \in
\cala$. Similarly, since $\varphi$ is multiplicative, we deduce
from (\ref{scal'}) and (\ref{add'}) that
\begin{eqnarray*}
\big(\varphi(\lambda b) -\lambda\varphi(b)\big)\varphi(a)
p&=&\varphi(\lambda b)\varphi(a) p-\lambda\varphi(b)\varphi(a)
p\\&=&\lambda\varphi(ba)p-\lambda\varphi(ba)p\\&=&0,
\end{eqnarray*}
and\begin{eqnarray*}\big(\varphi(b+c) -\varphi(b)
-\varphi(c)\big)\varphi(a)p&=&\varphi(b+c)\varphi(a)p
-\varphi(b)\varphi(a)p-\varphi(c)\varphi(a)p
\\&=&\varphi\big((b+c)a\big)p-\varphi(ba)p-\varphi(ca)p\\&=&0,\end{eqnarray*}for
all $a, b, c\in \cala$, and all $\lambda\in \mathbb{C}$. Therefore
$\varphi(a)(\calk)\subseteq \calk $. Since $\varphi$ is a
$*$-mapping, we conclude that $p\varphi(a)=\varphi(a)p$, for all
$a\in\cala$. Now define the mapping $\Phi:\cala \to B(\calh)$ by
$\Phi(a)=\varphi(a)p$. One can easily see that $\Phi$ is a
$*$-homomorphism, and so it is automatically continuous. Using
(\ref{commut'}), we obtain
\begin{eqnarray*}d(ab)&=&d(ab)p\\&=&\varphi(a)d(b)p+d(a)\varphi(b)p\\
&=&\varphi(a)pd(b)+d(a)\varphi(b)p\\&=&\Phi(a)d(b)+d(a)\Phi(b)
\qquad (a, b \in \cala).
\end{eqnarray*}
Hence $d$ is a $\Phi$-derivation. We deduce from Theorem 3.7 of
\cite{M-M} that $d$ is continuous.
\end{proof}

\section{$d$-continuity}

We start this section with the following definition.
\begin{definition}
Suppose that $\cala$ and $\calb$ are two normed algebras, and
$T,S:\cala\to \calb$ are two mappings.\\(i) The mapping $T$ is
\emph{left (resp. right) $S$-continuous }if $\displaystyle{\lim_{x
\to 0}}\big(T(x)S(b)\big)=0$, for all $b\in \cala$ (resp.
$\displaystyle{\lim_{x \to0}}\big(S(b)T(x)\big)=0$, for all $b\in
\cala$). If $T$ is both left and right $S$-continuous, then it is
simply called $S$-continuous.\\(ii) As for linear mappings, the
\emph{separating space} of a mapping $T$ is defined to be
$$\cals(T) :=\{b\in\calb\mid\,\, \exists \{a_n\}\subseteq \cala,\,\, a_n\to
0,\,\, T(a_n)\to b\}.$$ We notice that for a nonlinear mapping
$T$, this set is not necessarily a linear subspace and it may even
be empty. Recall that if $\cala$ and $\calb$ are complete spaces
and $T$ is linear, then the closed graph theorem implies that $T$
is continuous if and only if $\cals(T)=\{0\}$.
\end{definition}

If in Example \ref{example} we consider $\gamma$ and $\theta$ to
be $*$-mappings, then the continuity of $d$ together with
Corollary \ref{coro4.3} below imply that $\varphi$ and $\psi$ are
$d$-continuous mappings.

\begin{lemma}\label{separating}
Let $\cala$ be a normed algebra, let $\calx$ be a normed
$\cala$-bimodule, and let $d:\cala \to \calx$ be a
$(\varphi,\psi)$-derivation for two mappings $\varphi,\psi
:\cala\to \cala$. If $\varphi$ and $\psi$ are left and right
$d$-continuous mappings, respectively,
then$$a\big(\psi(bc)-\psi(b)\psi(c)\big)=0=\big(\varphi(bc)-\varphi(b)\varphi(c)\big)a\,,$$for
all $b,c\in\cala$ and all $a\in \cals(d)$.
\end{lemma}

\begin{proof}Suppose that $a\in
\cals(d)$. Then $a=\lim_{n\to \infty}d(a_n)$ for some sequence
$\{a_n\}$ converging to zero in $\cala$. By Lemma \ref{algebraic}
$(i)$, we have
\begin{eqnarray*}a\big(\psi(bc)-\psi(b)\psi(c)\big)
&=&\lim_{n\to\infty}d(a_n)\big(\psi(bc)-\psi(b)\psi(c)\big)\\
&=&\lim_{n\to\infty}\big(\varphi(a_nb)-\varphi(a_n)\varphi(b)\big)d(c)
\\&=&\lim_{n\to\infty}\varphi(a_n
b)d(c)\\&&-\lim_{n\to
\infty}\varphi(a_n)d(bc)+\lim_{n\to\infty}\varphi(a_n)d(b)\psi(c)\\&=&0
\qquad (b,c\in\cala),
\end{eqnarray*}
since $\varphi$ is left $d$-continuous.
\end{proof}

In the rest of this section, we assume that $\cala$ is a
$C^*$-subalgebra of $B(\calh)$, the $C^*$-algebra of all bounded
linear operators on a Hilbert space $\calh$. Also $\varphi, \psi$,
and $d$ are mappings from $\cala$ into $B(\calh)$. Removing the
assumption `linearity' and weakening the assumption `continuity'
on $\varphi$ and $\psi$, we extend the main result of \cite{M-M}
as follows.

\begin{theorem}\label{4.1.}
Let $\varphi$ be a $*$-mapping and let $d$ be a
$*$-$\varphi$-derivation. If $\varphi$ is left $d$-continuous,
then $d$ is continuous. Conversely, if $d$ is continuous then
$\varphi$ is left $d$-continuous.
\end{theorem}

\begin{proof}
Set
\begin{eqnarray*}
Y :=\{\varphi(bc)-\varphi(b)\varphi(c)\mid b,c\in\cala\}.
\end{eqnarray*} Since $\varphi$ is a $*$-mapping, $Y$ is a self-adjoint subset of $B(\calh)$. Put
$$\call_0
:=\bigcup\{u(h)\mid u\in Y, h\in \calh\}.$$ Let $\call$ be the
closed linear span of $\call_0$, and let $\calk :=\call^\perp$.
Then $\calh=\calk\oplus\call$, and by Lemma \ref{L perp=}, we
have\begin{eqnarray}\label{ker}\calk=\bigcap \{\ker(u)\mid u\in
Y\}.\end{eqnarray}Suppose that $p\in B(\calh)$ is the orthogonal
projection of $\calh$ onto $\calk$, then by (\ref{ker}), we
have\begin{eqnarray}\label{prod}\varphi(bc) p =
\varphi(b)\varphi(c)p,\end{eqnarray}for all $b,c\in \cala$. We
 claim that,
\begin{eqnarray}\label{commut}
p\varphi(a)=\varphi(a)p\,,\,\,\,pd(a)=d(a) p\qquad (a\in \cala).
\end{eqnarray}
To prove this, note that by Lemma \ref{algebraic}, and
(\ref{prod}), we have
\begin{eqnarray*}\big(\varphi(bc)
-\varphi(b)\varphi(c)\big)d(a)p=d(b)\big(\varphi(ca)-\varphi(c)\varphi(a)\big)p=0,
\end{eqnarray*}
for all $a, b, c \in\cala$. Thus $ud(a)p=0$, for all $a\in\cala$,
and all $u\in Y$. So $$ud(a)(\calk)=u d(a)p(\calh)=\{0\}.$$ This
means that $d(a)(\calk)$ is contained  in the kernel of each $u$
in $Y$. So by (\ref{ker}), $d(a)(\calk)\subseteq \calk$. Since
$d$ is a $*$-mapping, it follows that $pd(a)=d(a)p$. Similarly,
using (\ref{prod}), we get
\begin{eqnarray*}\big(\varphi(bc)
-\varphi(b)\varphi(c)\big)\varphi(a) p&=&\varphi(bc)\varphi(a)
p-\varphi(b)\varphi(c)\varphi(a)p\\&=&\varphi(bca) p-\varphi(bca)
p\\&=&0 \qquad (a, b, c\in \cala).
\end{eqnarray*}
 Therefore $\varphi(a)(\calk)=\varphi(a)p(\calh)\subseteq
\calk$. Since $\varphi$ is a $*$-mapping we conclude that
$p\varphi(a)=\varphi(a)p$, for all $a\in \cala$. Now define the
mappings $\Phi, D:\cala \to B(\calh)$ by $\Phi(a) :=\varphi(a)p$,
and $D(a) :=d(a)p$. We show that $\Phi$ is a multiplicative
$*$-mapping and $D$ is a $*$-$\Phi$-derivation. Clearly
(\ref{prod}) implies that $\Phi$ is multiplicative, and by
(\ref{commut}) we have
$$\big(\Phi(a)\big)^*=\big(\varphi(a) p\big)^*=p^*
\big(\varphi(a)\big)^*=p \varphi(a^*)=\Phi(a^*)\qquad(a\in
\cala).$$ Thus $\Phi$ is a multiplicative $*$-mapping. Now for
$a, b \in \cala$,
\begin{eqnarray*}
D(ab)&=&d(ab)p\\&=&\varphi(a)d(b)p+d(a)\varphi(b)p\\
&=&\varphi(a)pd(b)p+d(a)p\varphi(b)p\\&=&\Phi(a)D(b)+D(a)\Phi(b).
\end{eqnarray*}
Thus $D$ is a $*$-$\Phi$-derivation, and it is continuous by
Theorem \ref{multipl}. Now, we show that $\cals(d)=\{0\}$. Let $a
\in \cals(d)$, then there exists a sequence $\{a_n\}$ converging
to $0$ in $\cala$ such that $d(a_n)\to a$ as $n \to \infty$. Take
$h=k+\ell\in\calk\oplus \call=\calh$. By Lemma \ref{separating},
and by the fact that each $a\in \cals(d)$ is a bounded operator on
$\calh$, and that $\ell\in \call$ is in the closed linear span of
elements of the form
$\big(\varphi(bc)-\varphi(b)\varphi(c)\big)h$, where $b,
c\in\cala, h\in \calh$, we have $a(\ell)=0$. It follows from
continuity of $D$ that
$$a(k)=a\big(p(h)\big)=\lim_{n\to\infty}d(a_n)\big(p(h)\big)=\lim_{n\to
\infty}D(a_n)(h)=0.$$ Thus $a(h)=a(k)+a(\ell)=0$, and so $d$ is
continuous.

Conversely, let $d$ be a continuous $\varphi$-derivation. Then
$$\lim_{x\to 0}\varphi(x)d(b)=\lim_{x\to
0}d(xb)-\lim_{x\to0}d(x)\varphi(b)=0 \quad(b\in
\cala).$$\end{proof}

\begin{theorem}\label{4.2.}
Every $\varphi$-derivation $d$ is automatically continuous,
provided that $\varphi$ is $d$-continuous, and at least one of
$\varphi$ or $d$ is a $*$-mapping.
\end{theorem}

\begin{proof}
Let $d: \cala \to B(\calh)$ be a $\varphi$-derivation. Then
clearly $d^*$ is a $\varphi^*$-derivation. Set
$$\varphi_1 :=\frac{1}{2}(\varphi+\varphi^*),\quad \varphi_2 :=\frac{1}{2{\bf i}}(\varphi-\varphi^*
),\quad d_1 :=\frac{1}{2}(d+d^*),\quad d_2 :=\frac{1}{2{\bf
i}}(d-d^*).$$ Obviously these are $*$-mappings,
$\varphi=\varphi_1+ {\bf i}\varphi_2$, $d=d_1+ {\bf i}d_2$, and
$\varphi_k$ is $d_j$-continuous for $1 \leq k, j \leq 2$. A
straightforward calculation shows that if $\varphi$ is a
$*$-mapping then $\varphi_1=\varphi=\varphi_2$, and $d_1 , d_2$
are $\varphi$-derivations. Similarly, if $d$ is a $*$-mapping,
then $d_1=d=d_2$ and $d$ is a $\varphi_j$-derivation for $j=1, 2$.
Since $\varphi$ or $d$ is a $*$-mapping, then $\varphi_k$ is a
$*$-mapping and $d_j$ is a $*$-$\varphi_k$-derivation for $1 \leq
k, j \leq 2$. By Theorem \ref{4.1.}, the $d_j$'s are continuous,
and so $d=d_1+{\bf i}d_2$ is also continuous.
\end{proof}

\begin{corollary}\label{coro4.3}
Let $\varphi$ and $\psi$ be $*$-mappings and let $d$ be a
$*$-$(\varphi,\psi)$-derivation. Then $d$ is automatically
continuous if and only if $\varphi$ and $\psi$ are left and right
$d$-continuous, respectively.
\end{corollary}

\begin{proof}
Suppose that $\varphi$ and $\psi$ are left and right
$d$-continuous, respectively. Since $d, \varphi$, and $\psi$ are
$*$-mappings, then both $\varphi$ and $\psi$ are $d$-continuous.
Hence $\frac{\varphi+\psi}{2}$ is also $d$-continuous.
We have \begin{eqnarray*}2d(ab)&=&d(ab)+d^*(ab)\\
&=&\varphi(a)d(b)+d(a)\psi(b)+\big(\varphi(b^*)d(a^*)+d(b^*)\psi(a^*)\big)^*\\
&=&\varphi(a)d(b)+d(a)\psi(b)+\psi(a)d(b)+d(a)\varphi(b)\\
&=&(\varphi+\psi)(a)d(b)+d(a)(\varphi+\psi)(b) \qquad (a,
b\in\cala).
\end{eqnarray*}
Thus $d$ is a $*$-$\frac{\varphi+\psi}{2}$-derivation. It follows
from Theorem \ref{4.2.} that $d$ is continuous.

Conversely if $d$ is a continuous $(\varphi,\psi)$-derivation,
then
$$\lim_{x\to 0}\varphi(x)d(b)=\lim_{x\to 0}d(xb)-\lim_{x\to
0}d(x)\psi(b)=0 \quad(b\in\cala),$$and $$\lim_{x\to
0}d(b)\psi(x)=\lim_{x\to0}d(bx)-\lim_{x\to 0}\varphi(b)d(x)=0
\quad(b\in \cala).$$ So $\varphi$ and $\psi$ are left and right
$d$-continuous, respectively.
\end{proof}

\begin{lemma}
\label{phi(0)=0} Let $d$ be a $(\varphi,\psi)$-derivation. Then
there are two mappings $\Phi$ and $\Psi$ with $\Phi(0) = 0 =
\Psi(0)$ such that $d$ is a $(\Phi,\Psi)$-derivation.
\end{lemma}

\begin{proof}
Define $\Phi$ and $\Psi$ by
\begin{eqnarray*}\Phi(a):=\varphi(a)-\varphi(0) ,\\
\Psi(a):=\psi(a)-\psi(0), \end{eqnarray*} for all $a \in \cala$.
We
have\begin{eqnarray*}&&0=d(0)=d(a\cdot0)=\varphi(a)d(0)+d(a)\psi(0)=d(a)\psi(0)
\qquad (a \in \cala),\\&&0=d(0)=d(0\cdot
a)=\varphi(0)d(a)+d(0)\psi(a)=\varphi(0)d(a) \qquad (a \in
\cala).\end{eqnarray*} Thus\begin{eqnarray*}
\Phi(a)d(b)+d(a)\Psi(b)=\varphi(a)d(b)+d(a)\psi(b)= d(ab) \qquad
(a,b \in \cala).
\end{eqnarray*}Hence $d$ is a
$(\Phi,\Psi)$-derivation.\end{proof}

\begin{corollary}\label{coro4.4}If
$*$-mappings $\varphi$ and $\psi$ are continuous at zero, then
every $*$-$(\varphi,\psi)$-derivation $d$ is automatically
continuous.
\end{corollary}

\begin{proof}Apply
Lemma \ref{phi(0)=0} and Corollary \ref{coro4.3}.
\end{proof}

Clearly the assumption of Corollary \ref{coro4.4} comes true
whenever $\varphi$ and $\psi$ are linear and bounded. The next
theorem states that when we deal with a continuous
$*$-$(\varphi,\psi)$-derivation, we may assume that $\varphi$ and
$\psi$ have `at most' zero separating spaces.

\begin{theorem}
Let $\varphi$ and $\psi$ be $*$-mappings. If $d$ is a continuous
$*$-$(\varphi,\psi)$-derivation, then there are $*$-mappings
$\varphi^{'}$ and $ \psi^{'}$ from $\cala$ into $ B(\calh)$ with
`at most' zero separating spaces such that $d$ is a
$*$-$(\varphi^{'},\psi^{'})$-derivation.
\end{theorem}

\begin{proof}
Set\begin{eqnarray*}Y&:=&\{d(a)\mid
a\in \cala\},\\
\call_0& :=&\cup\{d(a)h\mid a\in \cala, h\in \calh
\}.\end{eqnarray*} Let $\call$  be the closed linear span of
$\call_{0}$ in $\calh$, and let $\calk :=\call^\perp$.  Suppose
that $p \in B(\calh)$ is the orthogonal projection of $\calh$ onto
$\call$. It follows from continuity of operators $d(a)$ that
$d(a)(\call)\subseteq \call \quad (a\in \cala)$, and so
$$pd(a)=d(a)p\qquad (a\in \cala),$$ and
$$p\varphi(a)=\varphi(a)p, \qquad p~\psi(a)=\psi(a)p \qquad
(a\in \cala).$$ For a typical element $\ell=d(b)h$ of $\call_0$,
we have
\begin{eqnarray*}
\varphi(a)\ell&=&\varphi(a)d(b)h\\&=&d(ab)h-d(a)\psi(b)h\in \call,
\end{eqnarray*}
thus $\varphi(a)(\call_0)\subseteq \call$. Hence
$\varphi(a)(\call)\subseteq \call$.\\ The same argument shows that
$\psi(a)(\call)\subseteq\call$. Now we define $\varphi^{'}$ and $
\psi^{'}$ from $\cala$ into $ B(\calh)$ by $\varphi^{'}(a)
:=\varphi(a)p \quad (a\in \cala)$, and $\psi^{'}(a) :=\psi(a)p
\quad (a\in \cala)$. Clearly $\varphi^{'}$ and $ \psi^{'}$ are
$*$-mappings. Furthermore, $d$ is a
$*$-$(\varphi^{'},\psi^{'})$-derivation. In fact for all
$a,b\in\cala$, and $h\in \calh$, we
have\begin{eqnarray*}d(ab)h&=&\varphi(a)d(b)h+d(a)\psi(b)h\\
&=&\varphi(a)p d(b)h+p
d(a)\psi(b)h\\&=&\varphi^{'}(a)d(b)h+d(a)\psi^{'}(b)h \,,
\end{eqnarray*}
since $p$ commutes with $d(a), \varphi(a)$, and $\psi(a)$. Suppose
that $\cals(\varphi^{'}) \neq \emptyset$ and
$a\in\cals(\varphi^{'})$. If
$a=\displaystyle{\lim_{n\to\infty}}\varphi^{'}(a_n)$ for some
sequence $\{a_n\}$ in $\cala$ converging to zero, then
\begin{eqnarray*}a(k)&=&\lim_{n\to
\infty}\varphi^{'}(a_n)k\\&=&\lim_{n\to\infty}\varphi(a_n)pk\\&=&\lim_{n\to
\infty}\varphi(a_n)0\\&=&0 \qquad (k\in \calk).\end{eqnarray*} Let
$\ell_0=d(b)h\in \call_0$, where $b \in \cala$ and $h \in \calh$.
Then
\begin{eqnarray*}
a(\ell_0)&=&\lim_{n\to\infty}\varphi^{'}(a_n)d(b)h\\&=&\lim_{n\to
\infty}d(a_nb)h-\lim_{n\to
\infty}d(a_n)\psi(b)h\\&=&0,\end{eqnarray*} and by continuity of
operator $a$, we get $a(\ell)=0 \quad (l\in\call)$, and hence
$a=0$. Similarly the separating space of $\psi^{'}$ is empty or
$\{0\}$.\end{proof}

\section{generalized $(\varphi,\psi)$-derivations}

In this section, we study the continuity of generalized
$(\varphi,\psi)$-derivations.

\begin{definition}
Suppose that $\calb$ is an algebra, $\cala$ is a subalgebra of
$\calb$, $\calx$ is a $\calb$-bimodule, $\varphi,\psi :\calb \to
\calb$ are mappings, and $d:\cala \to\calx$ is a
$(\varphi,\psi)$-derivation. A linear mapping $\delta:\cala \to
\calx$ is a \emph{generalized $(\varphi,\psi)$-derivation}
corresponding to $d$
if$$\delta(ab)=\varphi(a)d(b)+\delta(a)\psi(b)\quad (a,b\in
\cala).$$\end{definition}

\begin{proposition}\label{generalized}Suppose
that $\cala$ is a $C^*$-algebra acting on a Hilbert space
$\calh$. A generalized $*$-$(\varphi,\psi)$-derivation
$\delta:\cala\to B(\calh)$ corresponding to the
$(\varphi,\psi)$-derivation $d$ is automatically continuous
provided that $\varphi: \cala\to B(\calh)$ is a left
$d$-continuous $*$-mapping, and $\psi:\cala\to B(\calh)$ is both
a right $d$-continuous and a $\delta$-continuous
$*$-mapping.\end{proposition}

\begin{proof}
Suppose that $\{a_n\}$ is a sequence in $\cala$, and $a_n\to 0$ as
$n \to\infty$. By the Cohen factorization theorem, there exist a
sequence $\{b_n\}$ in $\cala$, and an element $c\in \cala$ such
that $a_n =cb_n\quad $, for all $n\in \mathbb{N}$, and $b_n\to 0$
as $n \to \infty$. By Corollary \ref{coro4.3}, $d$ is continuous,
so $d(a_n)=d(cb_n)\to 0$ as $n \to \infty$. A straightforward
computation shows that$$(\delta-d)(xy)=(\delta-d)(x)\psi(y)\qquad
(x, y\in
\cala).$$Thus\begin{eqnarray*}\delta(a_n)&=&(\delta-d)(a_n)+d(a_n)\\&=&(\delta-d)(c)\psi(b_n)+d(a_n),
\end{eqnarray*}
that converges to zero as $n \to \infty$, since $\psi$ is
right$(\delta-d)$-continuous.
\end{proof}

\begin{corollary}
Let $\cala$ be a $C^*$-algebra acting on a Hilbert space $\calh$.
Suppose that $\varphi, \psi:\cala\to B(\calh)$ are continuous at
zero. Then every generalized $*$-$(\varphi,\psi)$-derivation
$d:\cala\to B(\calh)$ is automatically continuous.
\end{corollary}

\begin{proof}
Using the same argument as in the proof of Lemma \ref{phi(0)=0},
we may assume that $\varphi(0)=0=\psi(0)$. Thus $\varphi$ and
$\psi$ are $S$-continuous for each mapping $S$. Now the result is
obtained from Theorem \ref{generalized}.
\end{proof}

\end{document}